\documentclass[12pt]{article}%
\usepackage{amsfonts}
\usepackage{amsmath,amssymb}
\usepackage[applemac]{inputenc}
\usepackage{amsmath,amssymb,fullpage}
\usepackage{color}
\usepackage{amsmath}
\usepackage{amssymb}
\usepackage{amsmath}
\usepackage{graphicx}
\usepackage{tikz}%
\providecommand{\U}[1]{\protect\rule{.1in}{.1in}}
\usetikzlibrary{arrows,shapes,automata,petri}

\newcommand{\dproof}{\noindent {Proof.} \quad}
\newcommand{\fproof}{\hfill $\square$ \bigskip}
\newtheorem{definition}{Definition}[section]
\newtheorem{example}{Example}[section]
\newtheorem{theorem}[definition]{Theorem}
\newtheorem{problem}[definition]{Problem}
\newtheorem{remark}[definition]{ \it Remark}

\newtheorem{lemma}[definition]{Lemma}

\numberwithin{equation}{section}

\def\1B{\text{1\!\!I}}

\def\H{\mathbb{H}}
\begin{document}

\title{SPDEs with space interactions and application to population modelling}
\author{Nacira Agram$^{1}$, Astrid Hilbert$^{1}$, Khouloud Makhlouf $^{2}$ and Bernt \O ksendal$^{3}$}
\date{30 June 2021}
\maketitle

\begin{abstract}
We consider optimal control of a new type of non-local stochastic partial differential
equations (SPDEs). The SPDEs have \emph{space interactions}, in the sense that the dynamics of the system at time $t$ and position in space $x$ also
depend on the space-mean of values at neighbouring points. This is a model
with many applications, e.g. to population growth studies and epidemiology.  We prove the existence and uniqueness of solutions of a class of SPDEs with space interactions, and we show that, under some conditions, the solutions are positive for all times if the initial values are. Sufficient and necessary maximum principles for the optimal control of
such systems are derived. Finally, we apply the results to study an optimal vaccine strategy problem
for an epidemic by modelling the population density as a space-mean stochastic reaction-diffusion equation.

\end{abstract}


\par
\footnotetext[1]{Department of Mathematics, Linnaeus University, V{\"a}xj{\"o}, Sweden. \\
Email: nacira.agram@lnu.se, astrid.hilbert@lnu.se.}
\par
\footnotetext[2]{Department of Mathematics, University of Biskra, Algeria. Email: khouloud.makhlouf@univ-biskra.dz.}
\footnotetext[3]{Department of Mathematics, University of Oslo, Norway. Email:  oksendal@math.uio.no.
}

\section{Introduction}

The purpose of this paper is to introduce a new type generalised stochastic heat equation with \emph{space interactions} as a model for population growth. By space interactions we mean that the dynamics of
the population density $Y(t,x)$ at a time $t$ and a point $x$ depends not only on its value and
derivatives at $x$, but also on its values in a neighbourhood of
$x$. For example, define $G$ to be a space-averaging operator of the form
\begin{equation}
G(x,\varphi)=\frac{1}{V(K_{r})}\int_{K_{r}}\varphi(x+y)dy;\quad
\varphi\in L^{2}(\mathbb{R}^n), \label{G}%
\end{equation}
where $V(\cdot)$ denotes Lebesgue volume and
\[
K_{r}=\{y\in\mathbb{R}^{n};|y|<r\}
\]
is the ball of radius $r>0$ in $\mathbb{R}^{n}$ centred at $0$. Then 
\begin{align*}
    \overline{Y}_G (t,x):= G(x,Y(t,\cdot))
\end{align*}
is the average value of $Y(t,x+ \cdot)$ in the ball $K_{r}$. \\
More generally, if we are given a nonnegative measure (weight) $\rho(dy)$ of total mass 1, then the $\rho$-weighted average of $Y$ at $x$ is defined by
\begin{align*}
    \overline{Y}_{\rho}(t,x):=\int_D Y(t,x+y)\rho(dy).
\end{align*}
\\
We believe that by allowing interactions between populations at different locations, we get a better model for population growth, including the modelling of epidemics. For example, we know that COVID-19 is spreading by close contact in space.

We illustrate the above by the following population growth model:
\begin{example}
With $G$ as in \eqref{G}, suppose the density $Y(t,x)$ of a population at the
time $t$ and the point $x$ satisfies the following space-interaction version of a
reaction-diffusion equation:
\begin{equation}%
\begin{cases}
dY(t,x) & =\Big(\frac{1}{2}\Delta Y(t,x)+\alpha\overline{Y}%
(t,x)-u(t,x)Y(t,x)\Big)dt+\beta{Y}(t,x)dB(t),\\
Y(0,x) & =\xi(x);\quad x\in D,\\
Y(t,x) & =\eta(t,x);\quad(t,x)\in(0,T)\times\partial D,
\end{cases}
\label{eq1.4*}%
\end{equation}
where $\alpha$ is a constant, $\xi,\eta$ are given bounded functions, $\overline{Y}(t,x)=G(x,Y(t,\cdot))$ and $B(t)=B(t,\omega); (t,\omega) \in [0,T] \times \Omega$ is a Brownian motion on a filtered probability space $(\Omega, \mathcal{F}, \{\mathcal{F}_t\}_{t \geq 0},P)$.\\
Here $u(t,x)$ is our control process, e.g. representing our harvesting or vaccine effort. 
\[
\begin{tikzpicture}
\coordinate (A1) at (4,0);
\coordinate (A2) at (5,1,);
\coordinate (A3) at (6.5,2);
\coordinate (A4) at (8.5,3);
\coordinate (A5) at (7,5);
\coordinate (A6) at (7,7);
\coordinate (A7) at (4,6);
\coordinate (A7a) at (3,5);
\coordinate (A8) at (2,4);
\coordinate (A9) at (0,2);
\draw[rounded corners=3.0mm] (A1) to [bend right] (A2) to [bend left] (A3) to
[bend right] (A4) to [bend right] (A5) to
[bend right] (A6) to [bend right] (A7) to
[bend left] (A8) to [bend right] (A9) to [bend right] (A1)
;
\draw[color=black, fill=lightgray] (5,3.5) circle [radius = 1.4];
\node at (5.5,3.5) {$x$};
\node at (5,5.8){$Y(t,x)$};
\node at (9,5) {$D$};
\end{tikzpicture}
\]

Then \eqref{eq1.4*} is a natural model for population growth in an environment
with space interactions. \newline

If $u(t,x)$ represents a vaccination effort rate at $(t,x)$, we define the total
expected utility $J_{0}(u)$ of the harvesting by an expression of the form
\begin{align*}
J_{0}(u)=\mathbb{E}\Big[\int_{D}\int_{0}^{T}U_{1}(u(t,x))dtdx+\int_{D}
U_{2}(Y(T,x))dx\Big],
\end{align*}
where $U_{1}$ and $U_{2}$ are given cost functions. The problem to find the
optimal vaccination rate $u^{\ast}$ is the following:

\begin{problem}
Find $u^{\ast}\in\mathcal{U}$ such that
\[
J_{0}(u^{\ast})=\inf_{u\in\mathcal{U}}J_{0}(u),
\]
where $\mathcal{U}$ is a given family of admissible controls.
\end{problem}
\end{example}

We will return to the example above after first discussing more general stochastic optimal control models 
with a system whose state $Y(t,x)$ at time $t$ and at the point $x$ satisfies an
SPDE with a non-local space-interaction dynamics of the following type:
\begin{equation}
\left\{
\begin{array}
[c]{ll}%
dY(t,x) & =A_{x}Y(t,x)dt+b(t,x,Y(t,x),Y(t,\cdot),u(t,x))dt\\
& +\sigma(t,x,Y(t,x),Y(t,\cdot),u(t,x))dB(t)
,\\
Y(0,x) & =\xi(x);\quad x\in D,\\
Y(t,x) & =\eta(t,x);\quad(t,x)\in(0,T)\times\partial D.
\end{array}
\right.  \label{eq1.1}%
\end{equation}
Here $dY(t,x)$ denotes the differential with respect to $t$ while $A_{x}$ is
the second order partial differential operator acting on $x$ of the form
\begin{equation}
A_{x}\phi(x)=\sum_{i,j=1}^{n}\alpha_{ij}(x)\frac{\partial^{2}\phi}{\partial
x_{i}\partial x_{j}}+\sum_{i=1}^{n}\beta_{i}(x)\frac{\partial\phi}{\partial
x_{i}};\quad\phi\in \mathcal{C}_{0}^{2}(%
%TCIMACRO{\U{211d} }%
%BeginExpansion
\mathbb{R}
%EndExpansion
^{n}).\label{opr}%
\end{equation}
Precise conditions on the coefficients will be given in the beginning of Section 4.2.

The domain $D$ is an open set in $\mathbb{R}^{n}$ with a
Lipschitz boundary $\partial D$ and closure $\overline{D}$. We extend $Y(t,x)$
to be a function on all of $[0,T]\times\mathbb{R}^{n}$ by setting
\[
Y(t,x)=0\text{ for }x\in\mathbb{R}^{n}\setminus\overline{D}.
\]

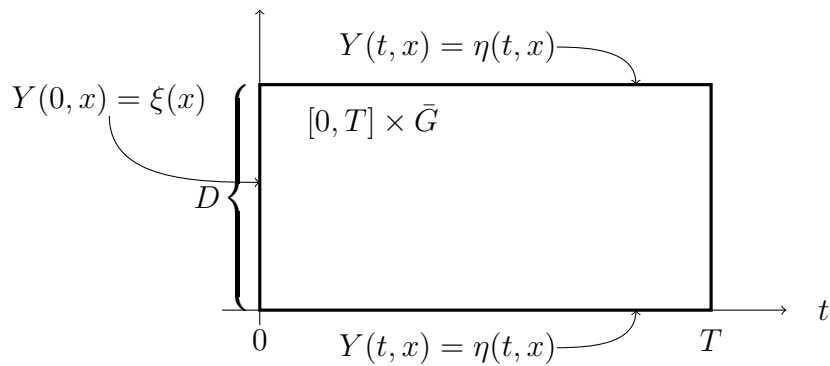
\[
\begin{tikzpicture}
\coordinate (LL) at (0,0);
\coordinate (UR) at (6,3);
\draw [very thick] (LL) rectangle (UR);
\draw [->] (-0.5,0) to (7,0);
\draw [->] (0,-0.2) to (0,4);
\node at (0,-0.4) {$0$};
\node at (6,-0.4) {$T$};
\node at (7.5,0) {$t$};
\node (Text1)[inner sep=0pt] at (-2,2.8) {$Y(0,x) = \xi(x)$};
\coordinate (dummy1) at (0,1.7);
\draw [->] (Text1.south) to [out=270,in=180] (dummy1.west);
\draw (-0.3,1.5) node[yscale=1.0]   {\rotatebox{270}{\makebox[3cm]{\upbracefill}}};
\node at (-0.7,1.5) {$D$};
\node (Text2)[inner sep=0pt] at (2.5,3.5) {$Y(t,x) = \eta(t,x)$};
\coordinate (dummy2) at (5,3.0);
\draw [->] (Text2.east) to [out=0,in=90] (dummy2.north);
\node at (1.5,2.5) {$[0,T] \times \bar D$};
\node (Text3)[inner sep=0pt] at (2.5,-0.5) {$Y(t,x) = \eta(t,x)$};
\coordinate (dummy3) at (5,0.0);
\draw [->] (Text3.east) to [out=0,in=270] (dummy3.south);
\end{tikzpicture}
\]

\begin{example}
In particular, the partial differential operator $A_{x}$ could be the
Laplacian $\Delta$. or more generally an operator of the $div-grad$-form
\[
A_{x}(\varphi)=div(\alpha(x)\nabla\varphi)(x);\quad\varphi\in \mathcal{C}^{2}(D),
\]
where $div$ denotes the divergence operator, $\nabla$ denotes the gradient
and
\[
\alpha(x)=[\alpha_{i,j}(x)]_{1\leq i,j\leq n}\in\mathbb{R}^{n\times n}%
\]
is a nonnegative definite matrix for each $x$. Equations of this type are of
interest because they represent important models in many situations, e.g. in
physics (e.g. fluid flow in random media, see e.g Holden \textit{et
al} \cite{HOUZ}), in epidemiology and in biology, e.g. in population growth
where $Y(t,x)$ represents the population density at $t,x$.
\end{example}

The new feature with this paper, is that we in addition to the operator
$A_{x}$ also allow a space dependence in the dynamics of the equation,
represented by the term $Y(t,\cdot)$ in \eqref{eq1.1}.

\begin{remark}
With the control $u$ given, sufficient conditions for the existence and uniqueness of the solution of the corresponding SPDEs with space interactions are known from general results on SPDEs. See e.g. Theorem 3.3  in Gawarecki \& Mandrekar \cite{GM}. However, very little seems to be known so far about the properties of such solutions. In Sections 2 and 3 we prove that the solution of a class of space-interaction reaction-diffusion equations can be obtained as a limit of an iterative procedure, and it is positive if the initial values are. See e.g. Theorem 2.1 and Theorem 3.3. This is an important confirmation that such equations are suitable models for population growth in general.
\end{remark}

There are two well-known approaches to solve stochastic control problems: The Bellman
dynamic programming method and the Pontryagin maximum principle.
Because of the space-mean dependence in our model, the system is not Markovian, and  it is not clear how to apply a dynamic programming approach. In stead we will use a stochastic version of the Pontryagin maximum principle, which involves a coupled system of a forward/backward SPDEs,

In the classical case when there is no interaction from neighbouring places, stochastic control of SPDEs has been studied widely in the literature, for example, we refer to Bensoussan \cite{b1}, \cite{b2}, \cite{b3},
\cite{b4}, Hu \& Peng \cite{HP}, Zhou \cite{Z}, \O ksendal \cite{o}, Fuhrman et a \cite{FHT} and \O ksendal \textit{et al } \cite{OPZ}, \cite{OSZ}, \cite{OS} and the references therein. \newline 
In the case of a control problem
for an SPDE with space-interaction dynamics we derive an adjoint process, which is a
backward SPDE with space-interaction dependence. For related singular stochastic control with space-interaction, we refer to Agram \textit{et al }\cite{AHO}.
\newline More details about the
theory of SPDE, we refer for example to Gawarecki \& Mandrekar \cite{GM}, Da Prato \& Zabczyk \cite{PZ},
Pardoux \cite{pard}, \cite{pardoux}, Hairer \cite{Hairer}, Pr\'ev\^ ot \& Roeckner
\cite{PR} and to Roeckner \& Zhang \cite{RZ}. \\
%\textcolor{red}{Most of what will be presented in this paper generalizes easily to jumps and to %multi-dimensions.}
Here is a summary of the content of this paper:

\begin{itemize}
\item In Section 2 we prove the existence and uniqueness of the solution of a class of space-interaction SPDEs, including the application studied in Section 5, and we give an iterative procedure for finding the solution (Theorem 2.1).
This result is new.
\item In Section 3 we use white noise theory to prove a useful positivity theorem for a class of SPDEs with space interactions (Theorem 3.1), and we prove that the solution is positive if the initial values are (Theorem 3.2).These results are also new and of independent interest.

\item Subsequently, in Section 4 we study the general optimization problem for such a
system. We derive both sufficient and necessary maximum principles for the
optimal control. See Theorem 4.5 and Theorem 4.6.

\item Finally, as an illustration of our results, in Section 5 we study an
example about optimal vaccination strategy for an epidemics modelled as an SPDE with space-interactions.
\end{itemize}

\section{Solutions of SPDEs with space interactions, \\
and positivity}

In this section we prove an existence and uniqueness result for solutions of SPDEs with space interactions. We are not aiming at proving this for the most general SPDE of this type, but we settle for a class of SPDEs which includes the application in Section 5. Thus, for simplicity we consider only the case when $A_x=L$ given by
\begin{align*}
L=\frac{1}{2}\Delta:=\frac{1}{2}\sum_{k=1}^{k=n}\frac{\partial^{2}}{\partial x_{k}^{2}},
\text{ and }D = \mathbb{R}^n,
\end{align*}
but it is clear that our method can also be applied to more general situations. \\
Fix $t>0,$ and let  $k\in\mathbb{N}_{0}=\left\{ 0,1,2,\ldots,\ldots\right\} ,$ $%
\alpha=\left( \alpha_{1},\alpha_{2},\ldots,\alpha_{m}\right) \in \mathbb{N}%
_{0}^{m}; m=1,2, ....$ \\
For functions $f \in \mathcal{C}_0^{\infty}(\mathbb{R}^n)$ (the family of functions in $\mathcal{C}(\mathbb{R}^n)$ with compact support), we define the Sobolev norm
\begin{equation*}
\left\vert f\right\vert _{k}=\sum_{|\alpha| \leq k}\big(\int_{\mathbb{R}^n} |\partial ^{\alpha} f(x)|^2 dx\big)^{\frac{1}{2}}; \alpha =\left( \alpha _{1},\alpha
_{2},\ldots ,\alpha _{n}\right) \in \mathbb{N}_{0}^{n},  \label{2b}
\end{equation*}%
and we define the Sobolev space $\H^k$ to be the closure of $\mathcal{C}_0^{\infty}(\mathbb{R}^n)$ in this norm. \\
Note that $\H^k$ is a Hilbert space for all $k$.\\

Also,  note that if $f \in \H^{k+2}$ then $Lf \in \H^k$, because
\begin{align}\label{2.1a}
|Lf|_k= \sum_{|\alpha| \leq k} \big(\int_{\mathbb{R}^n} |\partial ^{\alpha} Lf(x)|^2 dx\big)^{\frac12} 
\leq \frac12  \sum_{|\alpha| \leq k+2} \big(\int_{\mathbb{R}^n} |\partial ^{\alpha} f(x)|^2 dx\big)^{\frac12} =\frac12 |f|_{k+2}. 
\end{align}
%We let $\H$ denote the intersection (projective limit) of all the spaces $\H^{k}; k=0,1,...$.\\

Let $\mathcal{Y}_{k}^{\left( t\right) }$ denote the family of adapted random fields $%
Y\left( s,x\right) =Y\left( s,x,\omega\right) ,$ such that $||Y||_{k}^{(t)} < \infty$ where 
\begin{equation}
||Y||_{t,k}=\mathbb{E}\left[
\sup_{s\leq t}\left\{ \left\vert Y\left( s,.\right) \right\vert _{k}^{2}\right\} \right] ^{
{\frac12}},  \label{1b}
\end{equation}
and let $\mathcal{Y}^{(t)}$ be the intersection of all the spaces $\mathcal{Y}_{k}^{(t)}; k \in \mathbb{N}_0$, with the norm
\begin{align}
|| Y||_t^2:= \sum_{k=1}^{\infty} 2^{-k}||Y||_{t,k}^2.
\end{align}
In the following we let
$$\varphi \mapsto \overline{\varphi}(x)$$ 
be any averaging operator such that there exists a constant $C_1$ such that
\begin{align} \label{bar}
    |\overline{\varphi}|_{k} \leq C_1 |\varphi|_{k} \text{ for all } \varphi,k.
\end{align}
This holds, for example, if
$\overline{\varphi}(x)=\int \varphi(x+y) \rho(dy)$ for some measure $\rho$ of total mass 1.\\

We can now prove the following:
\begin{theorem}\label{theorem2.1}
Let $\xi \in \mathcal{Y}^{(T)}$ be deterministic and let $h:[0,T] \mapsto \mathbb{R}$ be bounded and deterministic. 
\begin{itemize}
\item[(i)]
Then there exists a unique solution $Y(t,x)\in \mathcal{Y}^{(T)}$ of the following SPDE with space interactions:
\begin{align*}
    Y(t,x) &= \xi(x)+\int_0^t LY(s,x)ds \nonumber\\ 
    &+ \int_0^t \overline{Y}(s,x) ds + \int_0^t h(s)Y(s,x)dB(s);\quad t \in [0,T].
\end{align*}
\item[(ii)] Moreover, the solution $Y(t,x)$ can be found by iteration, as follows:\\ 
Choose $Y_{0}\in \mathcal{Y}^{(T)}$ arbitrary deterministic and define inductively $Y_m$ to be the solution of
\begin{align}\label{recursive}
    Y_m(t,x) &= \xi(x)+\int_0^t LY_m(s,x)ds + \int_0^t \overline{Y}_{m-1}(s,x) ds \nonumber\\
    &+ \int_0^t h(s)Y_m(s,x)dB(s);\quad t \in [0,T]; m=1,2, ....
\end{align}
Then%
\begin{equation*}
Y_{m}\rightarrow Y \text{ in } \mathcal{Y}^{(T)} \text{ when } m \rightarrow \infty .
\end{equation*}

\end{itemize}
\end{theorem}

\dproof
(i): For $i=1,2$ choose $Z_i\in \mathcal{Y}^{(T)}$ and define 
$Y_i=Y^{Z_i}=: F(Z_i) $ to be the solution of the SPDE%
\begin{equation*}
Y^{Z_i}\left( t,x\right) =\xi \left( x\right) +\int_{0}^{t}LY^{Z_i}\left(
s,x\right) ds+\int_{0}^{t}\overline{Z_i}\left( s,x\right)
ds+\int_{0}^{t}Y^{Z_i}\left( s,\cdot\right) h\left( s\right) dB\left( s\right) .
\label{3b}
\end{equation*}
Note that here $Z_i$ (and hence $\overline{Z}_i)$, is given for each $i$. Therefore the existence and uniqueness of the solution  $Y^{Z_i}$ follows by the general existence and uniqueness theorems for solutions of SPDEs. e.g. as given in Theorem 3.3 in \cite{GM}.
Define%
\begin{align*}
&\widetilde{Y}=Y_{1}-Y_{2},\\
&\widetilde{Z}\left( t,x\right) =Z_{1}\left( t,x\right) -Z_{2}\left(
t,x\right) \\
&\widetilde{\overline{Z}}=\overline{Z}_{1}-\overline{Z}_{2}.
\label{5b}
\end{align*}%
Then%
\begin{equation*}
\widetilde{Y}\left( t,x\right) =\int_{0}^{t}L\widetilde{Y}\left( s,x\right)
ds+\int_{0}^{t}\widetilde{\overline{Z}}\left( s,x\right) ds+\int_{0}^{t}%
\widetilde{Y}\left( s,x\right) h\left( s\right) dB\left( s\right) .
\label{6b}
\end{equation*}%
Hence%
\begin{align}  
\left\vert \widetilde{Y}\left( s,.\right) \right\vert _{k}
& \leq
\int_{0}^{s}\left\vert L\widetilde{Y}\left( r,.\right) \right\vert
_{k}dr+\int_{0}^{s}\left\vert \widetilde{\overline{Z}}\left(
r,.\right) \right\vert _{k}dr \nonumber \\
&+\left\vert \int_0^s \widetilde{Y}(r,\cdot) h(r) dB(r)\right\vert_k \label{7b}
\end{align}%
By \eqref{2.1a} we have 
\begin{equation}
\left\vert L\widetilde{Y}\left( r,.\right) \right\vert _{k}\leq
\left\vert \widetilde{Y}\left( r,.\right) \right\vert _{k+2},
\label{8b}
\end{equation}%
and from \eqref{bar} we get
\begin{equation}
\left\vert \overline{Z}\right\vert _{k}\leq C_{1}\left\vert
Z\right\vert _{k }\text{ for all }k.  \label{9b}
\end{equation}%
Then by \eqref{7b}, \eqref{8b} and \eqref{9b}, we get%
\begin{eqnarray}
&&\mathbb{E}\left[ \sup_{s\leq t}\left\vert \widetilde{Y}\left( s,.\right)
\right\vert _{k}^{2}\right]   \label{10b} \\
&\leq &3\mathbb{E}\left[ \sup_{s\leq t}\left( \int_{0}^{s}\left\vert 
\widetilde{Y}\left( r,.\right) \right\vert _{k+2}dr\right) ^{2}\right] 
\notag +3C_{1}\mathbb{E}\left[ \sup_{s\leq t}\left( \int_{0}^{s}\left\vert 
\widetilde{Z}\left( r,.\right) \right\vert dr\right) ^{2}\right]   \notag \\
&&+3\mathbb{E}\left[ \sup_{s\leq t}\left\vert \int_{0}^{s} \widetilde{Y}\left( r,\cdot\right) h\left(
r\right) dB\left( r\right) \right\vert_k ^{2}\right] .  \notag
\end{eqnarray}%
By the Burkholder-Davis-Gundy inequality for Hilbert spaces (see e.g.  \cite{MR}), there exists a constant $C_2$ such that 
\begin{align}
\mathbb{E}\left[ \sup_{s\leq t}\left\vert \int_{0}^{s} \widetilde{Y}\left( r,\cdot\right) h\left(
r\right) dB\left( r\right) \right\vert_k ^{2}\right] 
&\leq C_{2}\mathbb{E}\left[ \int_{0}^{t}\left\vert \widetilde{Y}\left(
r,.\right) \right\vert _{k}^{2}h^2\left( r\right) dr\right]   \notag \\
&\leq C_{2}h_{0}^{2}t\mathbb{E}\left[ \sup_{s\leq t}\left\vert \widetilde{Y}%
\left( s,.\right) \right\vert _{k}^{2}\right] ;\text{ where } h_{0}^{2}=\sup_{s\in %
\left[ 0, T\right] }\left\vert h\left( s\right) \right\vert ^{2}.  \notag
\end{align}
Combining the above we get, if $0 \leq t \leq 1$,
\begin{align}
&\mathbb{E}\left[ \sup_{s\leq t}\left\vert \widetilde{Y}\left( s,.\right)
\right\vert _{k}^{2}\right]  \\
&\leq 3 t^{2}\mathbb{E}\left[ \sup_{r\leq t}\left\vert \widetilde{Y}%
\left( r,.\right) \right\vert _{k+2}^{2}\right] +3C_{1}t^{2}\mathbb{E}%
\left[ \sup_{r\leq t}\left\vert \widetilde{Z}\left( r,.\right) \right\vert
_{k}^{2}\right]  \\
&+3C_{2}h_{0}^{2}t\mathbb{E}\left[ \sup_{r\leq t}\left\vert \widetilde{Y}%
\left( r,.\right) \right\vert _{k}^{2}\right]. 
%\leq C_3 t \mathbb{E}\left[ \sup_{s\leq t}\left\vert \widetilde{Y}\left( s,.\right)
%\right\vert _{k+2}^{2}\right].
\end{align}%
In other words,
\begin{align} \label{2.14}
||\widetilde{Y}||_{t,k}^2 \leq 3t^2||\widetilde{Y}||_{t,k+2}^2 + 3C_1 t^2 ||\widetilde{Z}||_{t,k}^2 + 3C_2 h_0^2 t ||\widetilde{Y}||_{t,k}^2.
\end{align}
Note that 
\begin{align*}
\sum_{k=1}^{\infty} 2^{-k}||\widetilde{Y}||_{t,k+2}^2= \sum_{j=3}^{\infty}2^{-(j-2)}||\widetilde{Y}||_{t,j}\leq 4 \sum_{j=3}^{\infty}2^{-j}||\widetilde{Y}||_{t,j}
\leq 4 \sum_{k=1}^{\infty}2^{-k}||\widetilde{Y}||_{t,k} = 4 ||\widetilde{Y}||_t.
\end{align*}
Therefore, by multiplying the terms in \eqref{2.14} by $2^{-k}$ and summing over $k$, we get 
\begin{align*} 
||\widetilde{Y}||_{t}^2=\sum_{k=1}^{\infty} 2^{-k}||\widetilde{Y}||_{t,k}^2  \leq 12 t^2||\widetilde{Y}||_{t}^2 + 3C_1 t^2 ||\widetilde{Z}||_{t}^2 + 3C_2 h_0^2 t ||\widetilde{Y}||_{t}^2,
\end{align*}
or
\begin{align*}
(1-12 t^2 - 3C_2 h_0^2 t) ||\widetilde{Y}||_{t}^2  \leq 3C_1 t^2 || \widetilde{Z}||_t^2.
\end{align*}
Hence, if $t_{0}>0$ is chosen so small that%
\begin{equation*}
\frac{3C_1 t_0 ^2}{1-12 t_0^2 - 3C_2 h_0^2 t_0}< 1,
\end{equation*}%
we obtain that the map%
\begin{equation*}
Z\rightarrow Y^{Z}=F\left( Z\right) 
\end{equation*}%
is a contraction on $\mathcal{Y}^{\left( t_{0}\right) }.$ 
Therefore, by the Banach fixed point theorem there
exists a fixed point $\widehat{Y}$ of this map. Then $\widehat{Y}$ solves
the SPDE%
\begin{equation*}
\left\{ 
\begin{array}{l}
d\widehat{Y}\left( t,x\right) =L\widehat{Y}\left( t,x\right) dt+\overline{%
\widehat{Y}}\left( t,x\right) dt+\widehat{Y}\left( t,x\right) h\left(
t\right) dB\left( t\right) ; \quad t\in \left[ 0,t_{0}\right],  \\ 
\widehat{Y}\left( 0,x\right) =\xi(x); \quad x \in \mathbb{R}^n.
\end{array}%
\right.   \label{14b}
\end{equation*}
Uniqueness follows by a similar argument.

Since the constants do not depend on $t_{0},$ we can repeat the argument
starting from $t_{0}$ and hence by induction obtain a solution $Y(t,x)\in \mathcal{Y}^{(2t_0)}$. Repeating this argument we thus obtain a solution $Y \in \mathcal{Y}^{(T)}$.\\This proves part (i).\\

Part (ii) follows by the Banach fixed point theorem on the Banach space $\mathcal{Y}^{(T)}$.
\fproof

\section{The non-homogeneous stochastic heat equation and positivity}
In this section we will prove positivity of the solutions $Y(t,x)$ of SPDEs of the form
\begin{equation}
\left\{
\begin{array}
[c]{ll}
dY\left(  t,x\right)  &=LYdt+K\left(  t,x\right)  dt+h\left(  t\right)
Y\left(  t\right)  dB\left(  t\right), \nonumber\\
Y(0,x)&=\xi(x); \quad x \in \mathbb{R}^n,
\end{array}
\right. 
\end{equation}
where the function $\xi \in \mathcal{Y}^{(T)}$ is deterministic and positive, $h:[0,T] \mapsto \mathbb{R}$ is bounded and deterministic and $K(t,x)=K(t,x,\omega): [0,T] \times \mathbb{R}^n \times \Omega \mapsto \mathbb{R}$ is a given positive random field.\\
To motivate our method, we first recall the following basic results about the classical heat equation:\\
Let $L=\frac{1}{2}\triangle$ and consider the equation 
\begin{equation}
\left\{
\begin{array}
[c]{ll}%
dY\left(  t,x\right)  & =LYdt+K\left(  t,x\right)  dt,\\
Y\left(  0,x\right)  & =\xi(x);\quad x \in \mathbb{R}^n,
\end{array}
\right.  \label{1a}%
\end{equation}
where $\xi \in \mathcal{Y}^{(T)}$ and $K \in L^2([0,T] \times \mathbb{R}^n)$ are given deterministic functions.
Define the operator
$P_{t}:L^{2}\left(  \mathbb{R}^{n}\right)  \rightarrow
L^{2}\left(  \mathbb{R}^{n}\right)  $ by
\begin{equation}
P_{t}f\left(  x\right)  =\int_{\mathbb{R}^{n}}\left(  2\pi t \right)  ^{-\frac
{n}{2}}f\left(  y\right)  \exp\left(  -\frac{\left\vert x-y\right\vert ^{2}%
}{2t}\right)  dy,\label{2a}%
\end{equation}
then%
\begin{equation*}
\frac{d}{dt}P_{t}f=L\left(  P_{t}f\right),  \label{3a}
\end{equation*}
and if we define%
\begin{equation*}
Y\left(  t,x\right)  =P_{t}\xi\left(  x\right)  +\int_{0}^{t}P_{t-s}\left(
K\left(  s,.\right)  \right)  \left(  x\right)  ds, \label{4a}%
\end{equation*}
we get%
\begin{align*}
\frac{d}{dt}Y\left(  t,x\right)   &  =L\left(  P_{t}\xi\right)  \left(
x\right)  +P_{0}\left(  K\left(  t,.\right)  \right)  \left(  x\right)
+\int_{0}^{t}L\left(  P_{t-s}\left(  K\left(  s,.\right)  \right)  \right)
\left(  x\right)  ds\label{5a}\\
&  =LY\left(  t,x\right)  +K\left(  t,x\right)  .\nonumber
\end{align*}
Hence%
\[
Y\left(  t,x\right)  \text{ solves the heat equation \eqref{1a}}.
\]
Next, consider the case%
\begin{equation*}
dY\left(  t,x\right)  =LYdt+K\left(  t,x\right)  dt+\theta\left(  t\right)
Y\left(  t,x\right)  dt. \label{6a}%
\end{equation*}
Multiply the equation by%
\begin{equation*}
Z\left(  t\right)  =\exp\left(  -\int_{0}^{t}\theta\left(  s\right)
ds\right)  . \label{7a}%
\end{equation*}
Then the equation becomes%
\[
d\left(  Z\left(  t\right)  Y\left(  t,x\right)  \right)  =L\left(  Z\left(
t\right)  Y\left(  t,x\right)  \right)  dt+Z\left(  t\right)  K\left(
t,x\right)  dt.
\]
Hence, if we put
\begin{equation*}
\widehat{Y}=Z\left(  t\right)  Y\left(  t,x\right), \label{8a}%
\end{equation*}
then $\widehat{Y}$ solves the equation
\begin{equation*}
\left\{
\begin{array}
[c]{ll}%
d\widehat{Y}\left(  t,x\right)  & =L\widehat{Y}dt+Z\left(  t\right)K\left(  t,x\right)  dt,\\
\widehat{Y}\left(  0,x\right)  & =\xi(x),
\end{array}
\right.  \label{9a}%
\end{equation*}
and we are back to the previous case.\\
Finally, consider the SPDE%
\begin{equation}
dY( t,x)  =LYdt+K\left(  t,x\right)  dt+h\left(  t\right)
Y\left(  t\right)  dB\left(  t\right),\label{10a}%
\end{equation}
where $h$ is a given bounded deterministic function and $K(t,x)$ is stochastic and adapted, and $\mathbb{E}[\int_0^T\int_{\mathbb{R}^n} K^2(t,x) dt dx] < \infty$.
We handle this case by using white noise calculus on the Hida space $(\mathcal{S})^{*}$ of stochastic distributions: We introduce \emph{white noise} $W_t \in (\mathcal{S})^{*}$  defined by $$W_{t}=\frac{d}{dt}B(t),$$ 
and then we see that equation \eqref{10a} can be written
\begin{equation*}
\frac{d}{dt}Y\left(  t,x\right)  =LY+K\left(  t,x\right)  +Y\left(  t\right)
h\left(  t\right)  \diamond W_{t},\label{11a}%
\end{equation*}
where $\diamond$ denotes Wick multiplication. We refer to e.g. \cite{DOP} for more information about white noise calculus.
If we Wick-multiply this equation by
\begin{equation*}
Z_{t}:=\exp^{\diamond}\left(  -\int_{0}^{t}h\left(  s\right)  dB\left(
s\right)  \right),\label{12a}%
\end{equation*}
where in general $\exp^{\diamond}(\phi)=\sum_{n=0}^{\infty} \frac{1}{n!}\phi^{\diamond n}; \phi \in (\mathcal{S})^{*}$ is the Wick exponential,
we get%
\begin{equation}
Z_{t}\diamond\frac{d}{dt}Y\left(  t,x\right)  =L\left(  Y\diamond
Z_{t}\right)  +K\diamond Z_{t}+Y\left(  t\right)  h\left(  t\right)  \diamond
W_{t}\diamond Z_{t}.\label{13a}%
\end{equation}
Now
\begin{equation}
\frac{d}{dt}\left(  Z_{t}\diamond Y\right)  =Z_{t}\diamond\frac{d}{dt}Y\left(
t\right)  -Y\left(  t\right)  \diamond Z_{t}\diamond h\left(  t\right)
W_{t},\label{14a}%
\end{equation}
and hence \eqref{13a} can be written as
\[
\frac{d}{dt}
\underbrace{\left(  Z_{t}\diamond Y_{t}\right) }_{\text{ $
\widehat{Y}_t
$}}
 =L\underbrace{\left(  Z_{t}\diamond Y_{t}\right) }_{\text{ $
\widehat{Y}_t
$}} +K\left(  t,x\right)  \diamond Z_{t}.
\]
This has the same form as \eqref{1a}. Hence the solution $\widehat{Y}$ is
\[
\widehat{Y}\left(  t,x\right)  =P_{t}\xi\left(  x\right)  +\int_{0}%
^{t}P_{t-s}\left(  K\left(  s,.\right)  \right)  \left(  x\right)  \diamond
Z_{s}ds.
\]
Now we go back from $\widehat{Y}$ to $Y$ and get the solution
\begin{align}
Y\left(  t,x\right)   &  =\widehat{Y}\left(  t,x\right)  \diamond
\exp^{\diamond}\left(  \int_{0}^{t}h\left(  s\right)  dB\left(  s\right)
\right)  \nonumber\\
&  =P_{t}\xi\left(  x\right)  \diamond\exp^{\diamond}\left(  \int_{0}%
^{t}h\left(  s\right)  dB\left(  s\right)  \right)  \nonumber\\
&  +\int_{0}^{t}P_{t-s}\left(  K\left(  s,.\right)  \right)  \left(  x\right)
\diamond\exp^{\diamond}\left(  \int_{s}^{t}h\left( r\right)  dB\left(
r\right)  \right) ds.\label{15a}
\end{align}
Note that
\[
\exp^{\diamond}\left(  \int_{0}^{t}h\left(  s\right)  dB\left(  s\right)
\right)  =\exp\left(  \int_{0}^{t}h\left(  s\right)  dB\left(  s\right)
-\frac{1}{2}\int_0^t h^{2}\left(  s\right)  ds  \right)  >0.
\]
Recall the Gjessing-Benth lemma (see \cite{G}, \cite{BG} or Theorem 2.10.6 in \cite{HOUZ} or Proposition 13 in \cite{B}), which states that
\[
\phi\diamond\exp^{\diamond}\left(  \int_{0}^{t}h\left(  s\right)  dB\left(
s\right)  \right)  =\left(  \tau_{-h}\phi \right)  \exp^{\diamond}\left(  \int
_{0}^{t}h\left(  s\right)  dB\left(  s\right)  \right)  ,
\]
where, for $\phi: \Omega \mapsto \mathbb{R}$, we define $\tau_{-h}\phi(\omega)=\phi(\omega-h); \omega \in \Omega$ to be the shift operator on $\Omega$.

Using this in \eqref{15a} we conclude that if
\[
\xi\geq0\text{ and }K\geq0\text{ then }Y\geq0.
\]
We summarize what we have proved as follows:

\begin{theorem}\label{theorem3.1}
Assume that $\xi \in \mathcal{Y}^{(T)}$ is deterministic, $ \mathbb{E}[\int_0^T\int_{\mathbb{R}^n} K^2(t,x) dt dx] < \infty$ and let $h:[0,T] \mapsto [0,T]$ be bounded deterministic.
\begin{enumerate}
    \item 
Then the unique solution $Y(t,x)\in \mathcal{Y}^{(T)}$ of the
non-homogeneous SPDE
\begin{align*}
dY\left(  t,x\right)  &=LYdt+K\left(  t,x\right)  dt+h\left(  t\right)
Y\left(  t\right)  dB\left(  t\right), \nonumber\\
Y(0,x)&=\xi(x); \quad x \in \mathbb{R}^n
\end{align*}
is given by
\begin{align*}
Y\left(  t,x\right)   
&  =(\tau_{-h}P_{t}\xi)(x)\exp^{\diamond}\left(  \int_{0}%
^{t}h\left(  s\right)  dB\left(  s\right)  \right)  \nonumber\\
&  +\int_{0}^{t}(\tau_{-h}P_{t-s}(K(s,.))(x)
\exp^{\diamond}(\int_{s}^{t}h(r)dB(r)) ds,
\end{align*}
where $\exp^{\diamond}(\int_s^t h(r) dB(r))=\exp(\int_s^th(r) dB(r)-\frac{1}{2}\int_s^t h^2(r)dr);\quad 0\leq s \leq t \leq T.$
\item
In particular, if $\xi(x) \geq 0$ and $K(t,x) \geq 0$ for all $(t,x) \in [0,T] \times \mathbb{R}^n$, then $Y(t,x) \geq 0$ for all $(t,x)\in [0,T] \times \mathbb{R}^n$.
\end{enumerate}
\end{theorem}
Combining this with Theorem \ref{theorem2.1} we get
\begin{theorem}{(Positivity)}\label{theorem3.2}
Assume that $\xi \in \mathcal{Y}^{(T)}$ is deterministic and let $h:[0,T] \mapsto \mathbb{R}$ be bounded and deterministic. Let $Y(t,x) \in \mathcal{Y}^{(T)}$ be the unique solution of the following SPDE with space interactions:
\begin{align}\label{18a}
    Y(t,x) &= \xi(x)+\int_0^t LY(s,x)ds + \int_0^t \overline{Y}(s,x) ds + \int_0^t h(s)Y(s,x)dB(s); \quad t \in [0,T],
\end{align} given by Theorem \ref{theorem3.1}.\\
Then if $\xi(x) \geq 0$ for all $x \in \mathbb{R}^{n}$, we have $Y(t,x) \geq 0$ for all $(t,x) \in [0,T] \times \mathbb{R}^{n}$.
\end{theorem}

\dproof
By Theorem \ref{theorem2.1} we know that the solution of  \eqref{18a} can be obtained as the limit when $m \rightarrow \infty$ of the sequence $Y_m(t,x)$ defined recursively by
the equation \eqref{recursive}. Then by Theorem \ref{theorem3.1}, part 2, we know that $Y_m(t,x) \geq 0$ for all $t,x,m$. We conclude that $Y(t,x) \geq 0$ for all $(t,x)$.
\fproof

\begin{remark} \label{remark3.3}
The results from this and the previous section can be extended to equations of the form
\begin{align*}\label{19a}
    dY(t,x) &= \big[LY(t,x)+\gamma(t,x)Y(t,x)\big]dt+\overline{Y}(t,x) dt+ h(t)dB(t);\quad  t \in [0,T],
\end{align*}
for a given adapted process $\gamma \in \mathcal{Y}^{(T)}$.To see this we apply the arguments above with the operator $L$ replaced by the operator $\widehat{L}$ defined by $\widehat{L}\varphi=L\varphi +\gamma \varphi; \varphi \in \mathcal{C}^{\infty}(\mathbb{R}^n)$. We omit the details. 
\end{remark}

\section{The optimization problem}

In general, if $\mathcal{X},\mathcal{Y}$ are two Banach spaces and
$F:\mathcal{X}\rightarrow\mathcal{Y}$ if Fr\'echet differentiable at $x\in\mathcal{X}$, then
we let $\nabla_{x}F$ denote the Fr\'echet derivative of $F$ at $x$. It is a linear operator from $\mathcal{X}$ to $\mathcal{Y}$ and the action of $\nabla_x F$ to $h \in \mathcal{X}$ is denoted by $\nabla_x F(h)=\langle \nabla_x F,h\rangle \in \mathcal{Y}$. Recall that 
if $F$ is Fr\'echet differentiable at $x$ with Fr\'echet derivative
$\nabla_{x}F$, then $F$ has a directional derivative 
\begin{align}
D_{x}F(h):= \lim_{\epsilon \rightarrow 0} \frac{1}{\epsilon}(F(x+\epsilon h) -F(x))
\end{align}
in all directions $h \in \mathcal{X}$ and
\begin{align}
D_x F(h) =\nabla_x F(h)=\langle \nabla_x F,h \rangle.
\end{align}
\noindent In particular, note that if $F$ is a linear operator, then
$\nabla_{x}F=F$ for all $x$.

\subsection{The Hamiltonian and the adjoint BSPDE}

We now give a general formulation of the problem we consider.\newline Let
$A_x$ be a linear second order partial differential operator given by 
\begin{equation}
A_{x}\phi(x)=\sum_{i,j=1}^{n}\alpha_{ij}(x)\frac{\partial^{2}\phi}{\partial
x_{i}\partial x_{j}}+\sum_{i=1}^{n}\beta_{i}(x)\frac{\partial\phi}{\partial
x_{i}};\quad\phi\in \mathcal{C}_{0}^{2}(%
%TCIMACRO{\U{211d} }%
%BeginExpansion
\mathbb{R}
%EndExpansion
^{n}).\label{opr}%
\end{equation}

Let $T>0$ and assume that the state $Y(t,x)$ at time $t\in\lbrack0,T]$ and at
the point $x\in\overline{D}:=D\cup\partial D$ satisfies the following non-local 
quasilinear stochastic heat equation:
\begin{equation}
\left\{
\begin{array}
[c]{ll}%
dY(t,x) & =A_{x}Y(t,x)dt+b(t,x,Y(t,x),Y(t,\cdot),u(t,x))dt\\
& +\sigma(t,x,Y(t,x),Y(t,\cdot),u(t,x))dB(t)
,\\
Y(0,x) & =\xi(x);\quad x\in{D},\\
Y(t,x) & =\eta(t,x);\quad(t,x)\in(0,T)\times\partial D.
\end{array}
\right.  \label{eq2.1}%
\end{equation}
%Let $L(\mathbb{R}^{n})$ denote the set of real measurable functions on$\mathbb{R}^{n}$. \\
\emph{
\textbf{Assumptions on} ($\alpha,\beta, b, \sigma$)
\begin{itemize}
  \item [(a)] $(\alpha_{ij}(x))_{1\leq i,j\leq n}$ is a given symmetric nonnegative definite $n\times n$ matrix with eigenvalues bounded away from
0 and with entries $\alpha_{ij}(x)\in \mathcal{C}^{4}(D)\cap \mathcal{C}(\overline
{D})$ for all $i,j=1,2,...,n.$
 \item [(b)] $\beta_{i}(x)\in \mathcal{C}^{3}(D)\cap \mathcal{C}(\overline{D})$ for all $i=1,2,...,n$.
 \item [(c)] The functions $b$ and $\sigma$ are $\mathbb{F}$-adapted, $\mathcal{C}^{2}$ with respect to $y$ and $u$ and admit uniformly bounded derivatives.
 \item [(d)]For each $t,x,y,u$ the functions
\begin{align*}
\varphi &  \mapsto b(t,x,y,\varphi,u):[0,T]\times D\times%
%TCIMACRO{\U{211d} }%
%BeginExpansion
\mathbb{R}
%EndExpansion
\times L(\mathbb{R}^{n})\times U\rightarrow%
%TCIMACRO{\U{211d} }%
%BeginExpansion
\mathbb{R}
%EndExpansion
,\\
\varphi &  \mapsto\sigma(t,x,y,\varphi,u):[0,T]\times D\times%
%TCIMACRO{\U{211d} }%
%BeginExpansion
\mathbb{R}
%EndExpansion
\times L(\mathbb{R}^{n})\times U\rightarrow%
%TCIMACRO{\U{211d} }%
%BeginExpansion
\mathbb{R},
\end{align*}
are $\mathcal{C}^{1}$ functionals on $L^{2}(D)=L^{2}(D,m)$, where $dm(x)=dx$ is the
Lebesgue measure on $\mathbb{R}^{n}$. 
\end{itemize}
}
We call the equation \eqref{eq2.1} a
\emph{stochastic partial differential equation with space-interactions.}
\newline In general, the formal adjoint $A^{\ast}$ of an operator $A$ is defined by the identity
\begin{equation*}
(A\phi,\psi)=(\phi,A^{\ast}\psi),\,\,\,\text{ for all }\phi,\psi\in
\mathcal{C}_{0}^{2}(D),\label{int by part}%
\end{equation*}
where $(\phi_{1},\phi
_{2}):=\langle\phi_{1},\phi_{2}\rangle_{L^{2}(D)}=%
%TCIMACRO{\dint _{\mathbb{R}}}%
%BeginExpansion
{\displaystyle\int_{D}}
%EndExpansion
\phi_{1}(x)\phi_{2}(x)dx$ is the inner product in $L^{2}(D)$ and $\mathcal{C}_{0}^{2}(D)$ is the set of twice differentiable functions with compact support in $D$. In our
case we have
\[
A_{x}^{\ast}\phi(x)=\sum_{i,j=1}^{n}\frac{\partial^{2}}{\partial x_{i}\partial
x_{j}}(\alpha_{ij}(x)\phi(x))-\sum_{i=1}^{n}\frac{\partial}{\partial x_{i}
}(\beta_{i}(x)\phi(x));\quad\phi\in \mathcal{C}^{2}(D).
\]
We interpret $Y$ as a weak (variational) solution to \eqref{eq2.1}, in the
sense that
\begin{align*}
\langle Y(t),\phi\rangle_{L^{2}(D)} &  =\langle\xi(x),\phi\rangle_{L^{2}
(D)}+\int_{0}^{t}\langle Y(s),A_{x}^{\ast}\phi\rangle_{L^{2}(D)} ds\\
&  +\int_{0}^{t}\langle b(s,Y(s)),\phi\rangle_{L^{2}(D)}ds+\int_{0}^{t}%
\langle\sigma(s,Y(s)),\phi\rangle_{L^{2}(D)}dB(s); \phi \in \mathcal{C}_0^2(D).
\end{align*}
For simplicity, in the above equation we
have not written all the arguments of $b,\sigma$. We can in the usual
way apply the It\^ o formula to solutions of such SPDEs.\newline The process $u(t,x)=u(t,x,\omega)$ is
our control process, assumed to have values in a given convex set $U\subset%
%TCIMACRO{\U{211d} }%
%BeginExpansion
\mathbb{R}
%EndExpansion
^{k}$. We assume that $u(t,x)$ is $\mathbb{F}$-predictable for all
$(t,x)\in(0,T)\times D$. We call the control process $u(t,x)$ admissible if
the corresponding SPDE with space-mean dynamics (\ref{eq2.1}) has a unique
 solution $Y\in \mathcal{Y}^{(T)}$. The set of admissible controls is denoted by $\mathcal{U}$.\newline The
performance functional (cost) associated to the control $u$ is assumed to have the
form
\begin{equation}\label{jcost}
J(u)=\mathbb{E}\Big[\int_{0}^{T}\int_{D}f(t,x,Y(t,x),Y(t,\cdot
),u(t,x))dxdt+\int_{D}g(x,Y(T,x),Y(T,\cdot))dx\Big];\quad u\in\mathcal{U}.
\end{equation}
\textbf{Assumptions on} ($f,g$)
\begin{itemize}
    \item [(e)] The function $f(t,x,y,\varphi,u)$ is $\mathbb{F}$-adapted and the function $g(x,y,\varphi)$ is $\mathcal{F}_T$-measurable. They are assumed to be bounded,  $\mathcal{C}^2$  with respect to $y,\varphi,u$, with uniformly bounded derivatives.
\end{itemize}
%for each $t,x,y,u$ the functions $\varphi\mapsto f(t,x,y,\varphi
%,u):[0,T]\times D\times%
%TCIMACRO{\U{211d} }%
%BeginExpansion
%\mathbb{R}
%EndExpansion
%\times L(\mathbb{R}^{n})\times U\rightarrow%
%TCIMACRO{\U{211d} }%
%BeginExpansion
%\mathbb{R}
%EndExpansion
%,$ and $\varphi\mapsto g(x,y,\varphi):D\times%
%TCIMACRO{\U{211d} }%
%BeginExpansion
%mathbb{R}
%EndExpansion
%\times L(\mathbb{R}^{n})\rightarrow%
%TCIMACRO{\U{211d} }%
%BeginExpansion
%\mathbb{R}
%EndExpansion
%,$ are $C^{1}$ functionals on $L^{2}(D)$.\newline
%\newline 
The general problem
we consider in this paper is the following:

\begin{problem}
Find $\widehat{u}\in\mathcal{U}$ such that
\begin{equation}
J(\widehat{u})=\inf_{u\in\mathcal{U}}J(u). \label{eq2.4}%
\end{equation}

\end{problem}

To study this problem we define the associated \emph{Hamiltonian} $H:\left[
0,T\right]  \times D\times\mathbb{R}\times L(\mathbb{R}^{n})\times
U\times\mathbb{R}\times\mathbb{R}\times\Omega\rightarrow
\mathbb{R}$ by

\begin{align}
H(t,x,y,\varphi,u,p,q)  &  :=H(t,x,y,\varphi,u,p,q,\omega)=f(t,x,y,\varphi
,u)+b(t,x,y,\varphi,u)p\nonumber\\
&  \quad\quad+\sigma(t,x,y,\varphi,u)q. \label{h}%
\end{align}
%We assume that $H,f,b,\sigma,\gamma$ and $g$ are continuously differentiable
%($C^{1}$) and they admit Fr\'echet derivatives with respect to $y$ and bounded
%Fr\'echet derivatives with respect to $\varphi$ and $u.$ \newline 
In general, if
$h:L^{2}(D)\mapsto L^{2}(D)$ is Fr\'echet differentiable map, then its Fr\'echet
derivative (gradient) at $\varphi\in L^{2}(D)$ denoted by $\nabla_{\varphi}h=\nabla h$ is a bounded linear map on the Hilbert space $L^2(D)$, and by the Riesz representation theorem we can represent it by a function $\nabla h(x,y) \in L^2(D \times D)$.
We denote the action of $\nabla h$ on a function $\psi\in L^{2}(D)$ by
$\left\langle \nabla h,\psi\right\rangle $.
\vskip 0.3cm
Hence
\begin{equation}
\left\langle \nabla h,\psi\right\rangle (x):=\int_{D}\nabla h(x,y)\psi(y)dy;\quad\text{ for all }\psi\in L^{2}(D). \label{eq2.6a}%
\end{equation}

\begin{remark}
 \begin{itemize}
 \item Note in particular that if $h:L^{2}(D)\mapsto L^{2}(D)$ is linear, then

$$\nabla h(x,y)=h(x,y)$$.

 \item Also note that from \eqref{eq2.6a} it follows by the Fubini theorem that
\begin{align}
&  \int_{D}\left\langle \nabla h,\psi\right\rangle (x)dx=\int_{D}%
\int_{D}\nabla h(x,y)\psi(y)dydx=\int_{D}\int_{D}%
\nabla h(y,x)\psi(x)dxdy\nonumber\\
&  =\int_{D}\left(  \int_{D}\nabla h(y,x)dy\right)
\psi(x)dx=\int_{D}\overline{\nabla}h(x)\psi(x)dx,\nonumber
\end{align}
where
\begin{equation}
\overline{\nabla} h(x):=\int_{D}\nabla h(y,x)dy. \label{eq2.8a}%
\end{equation}

\end{itemize}
\end{remark}

\begin{example}
\label{exp}

\begin{description}
\item[a)] Assume that $h:L^{2}(D)\mapsto L^{2}(D)$ is given by
\begin{equation}
h(\varphi)=\left\langle h,\varphi\right\rangle (x)=G(x,\varphi(\cdot
))=\frac{1}{V(K_{r})}\int_{K_{r}}\varphi(x+y)dy. \label{eq2.6}%
\end{equation}

\end{description}

Then
\[
\left\langle \nabla_{\varphi}h,\psi\right\rangle (x)=\left\langle
h,\psi\right\rangle (x)=\frac{1}{V(K_{r})}\int_{K_{r}}\psi(x+y)dy.
\]

Therefore $\nabla h(x,y)$ is given by the identity
\[
\int_{D}\nabla h(x,y)\psi(y)dy=\frac{1}{V(K_{r})}%
\int_{K_{r}}\psi(x+y)dy;\quad\psi\in L^{2}(D).
\]
Substituting $z=x+y$ this can be written
\[
\int_{D}\nabla_{\varphi}^{\ast}h(x,y)\psi(y)dy=\frac{1}{V(K_{r})}%
\int_{x+K_{r}}\psi(z)dz=\int_{D}\frac{\mathbf{1}_{x+K_{r}}%
(y)}{V(K_{r})}\psi(y)dy.
\]
Since this is required to hold for all $\psi$, we conclude the following:

\begin{description}
\item[b)] \label{lem}Suppose that $h$ is given by \eqref{eq2.6}. Then
\[
\nabla h(x,y)=\frac{\mathbf{1}_{x+K_{r}}(y)}{V(K_{r})},
\]
and
\begin{align*}
\overline{\nabla}_{\varphi} h(x)  &  =\int_{D}\nabla_{\varphi}h(y,x)dy=\frac{1}{V(K_{r})}\int_{D}\mathbf{1}_{y+K_{r}}%
(x)dy=\frac{1}{V(K_{r})}\int_{D}\mathbf{1}_{x-K_{r}}%
(y)dy\\
&  =\frac{V((x-K_{r})\cap D)}{V(K_{r})}=\frac{V((x+K_{r})\cap
D)}{V(K_{r})},\nonumber
\end{align*}
since $K_{r}=-K_{r}$.
\end{description}
\end{example}
We associate to the Hamiltonian the following backward SPDE 
\begin{equation} 
dp(t,x) =-\left[  A_{x}^{\ast}p(t,x)+\tfrac{\partial H}{\partial
y}(t,x)+\overline{\nabla} H(t,x)\right]dt+q(t,x)dB(t), \label{p}%
\end{equation}
with boundary/terminal values
\begin{align}
\left\{
\begin{array}
[c]{ll}%
  p(T,x)&=\tfrac{\partial g}{\partial y}(x)+\overline{\nabla} g(x);\quad x\in{D},\\
     p(t,x)&=0;\quad (t,x)\in(0,T)\times\partial D,     \end{array}
     \right.\label{p_T}
\end{align}
where we have used the simplified notation
\[
H(t,x)=H(t,x,y,\varphi,u,p,q)|_{y=Y(t,x),\varphi=Y(t,\cdot
),u=u(t,x),p=p(t,x),q=q(t,x)},
\]
and similarly we have used the notation $g(x)$ for $g(x,Y(T,x),Y(T,\cdot))$. Here $A_{x}^{\ast}$ denotes the adjoint of the operator $A_x$.

\begin{remark}
Note that in differential form the BSPDE \eqref{p} can be written
\begin{align*}
dp\left( t,x\right)  & =-\Big[ \sum\limits_{i,j=1}^{n}\alpha
_{ij}(x)  \frac{\partial ^{2}}{\partial x_{i}\partial x_{j}}
p( t,x) \\
&+\sum\limits_{i=1}^{n} \Big(-\beta_i(x)+2\sum\limits_{j=1}^{n}\frac{\partial}{\partial x_j}\alpha_{ij}(x)\Big) \frac{\partial}{\partial x_i}p(t,x) \\
& \Big(-\sum\limits_{i=1}^{n} \frac{\partial}{\partial x_i}\beta_i(x)+\sum\limits_{i,j=1}^{n} \frac{\partial ^{2}}{\partial x_{i}\partial x_{j}} \alpha_{ij}(x)\Big)p( t,x)+\Big( (\frac{\partial }{\partial y} + \overline{\nabla} )b(t,x)\Big)p(t,x) \\
&+\Big((\frac{\partial }{\partial y} 
+ \overline{\nabla} )\sigma(t,x) \Big)q(t,x) +\Big(\frac{\partial }{\partial y} + \overline{\nabla} \Big)f(t,x) \Big] dt + q(t,x) dB(t).
%\\p(T,x) &=\tfrac{\partial g}{\partial y}(x)+\overline{\nabla}_{\varphi}^{\ast}g(x) ; x \in D,\\
%p(t,x)&= 0; (t,x) \in (0,T) \times \partial D.
\end{align*}
To the best of our knowledge, the existence and uniqueness of a solution of \eqref{p}-\eqref{p_T} is not known in general. However, note that (for given $u$) the equation \eqref{p}, regarded as a BSPDE in the unknown $\mathcal{Y}^{(T)}\times \mathcal{Y}^{(T)}$-valued processes $(p,q)$, is linear. Therefore, in view of our general assumptions (a)-(e) above,  the existence and uniqueness of solution follows from e.g.Theorem 2.1 in \cite{DTZ}, provided that the terms $\overline{\nabla} b(t,x),\overline{\nabla} \sigma(t,x)$ and $\overline{\nabla} f(t,x)$ satisfy condition $(F_m)$ in \cite{DTZ}. To this end, it suffices that 
$b, \sigma$ and $f$ depend linearly on $\varphi$ and in a space-averaging manner, as in the example with $h$ in  \eqref{eq2.6} above. In particular, this holds in the application studied in Section 5.
\end{remark}
\begin{remark} Here, as in Sections 2 and 3,  we are primarily interested in strong solutions $(p,q) \in \mathcal{Y}^{(T)} \times \mathcal{Y}^{(T)}$, but weak solutions are also of interest. A pair $(p,q)$ of random fields is said to be a weak solution to the backward SPDE \eqref{p}-\eqref{p_T} if, for all  $\phi\in \mathcal{C}_{0}^{2}(D),$
\begin{align*}
\langle p(t,.),\phi(.)\rangle -\langle p(T,.),\phi(.)\rangle &=\int_{t}^{T}\langle A_{x}^{\ast}p(s,.),\phi(.)\rangle ds +\int_{t}^{T}\langle \tfrac{\partial H}{\partial
y}(t,.)+\overline{\nabla} H(t,.),\phi(.)\rangle ds\\
& -\int_{t}^{T}%
\langle q(s,.),\phi(.)\rangle dB(s);\quad  \text{a.s. for each t}  \in [0,T].
\end{align*}
Hence, we observe that $p$ admits the following mild representation
\small
$$
p(t,x) =P_{T-t}\Big(p(T,x)\Big)+\int_t^TP_{s-t}\Big(\tfrac{\partial H}{\partial
y}(t,x)+\overline{\nabla} H(t,x)\Big)ds-\int_t^TP_{s-t}\Big(q(s,x)\Big)dB(s);\quad 0\leq t \leq T,
$$
where $P_t$ denotes the semigroup of the operator $A^{*}$.
\end{remark}
%\textcolor{blue}{We can rewrite the BSPDE \eqref{p}-\eqref{p_T} as follows}
%\[
%\left\{ 
%\begin{array}{ll}
%dp\left( t,x\right)  & =-\Big[ \sum\limits_{i,j=1}^{n}\overline{\alpha }%
%^{ij}\left( t,x\right) \frac{\partial ^{2}}{\partial x_{i}\partial x_{j}}%
%p( t,x) -\sum\limits_{i=1}^{n}\overline{\beta }^{i}( t,x) p(t,x) +\overline{b}( t,x) p( t,x) \\ 
%& +c\left( t,x\right) q\left( t,x\right) +\overline{f}\left(t,x\right) \Big] dt+q\left( t,x\right) dB\left( t\right)  \\ 
%p\left( T,x\right)  & =\overline{g}\left( T,x\right) \end{array}%
%\right. 
%\]
%where 
%\[
%\overline{g}\left( T,x\right)=\tfrac{\partial g}{\partial y}(x)+\overline{\nabla}_{\varphi}^{\ast}g(x), \overline{f}\left(
%t,x\right)=\tfrac{\partial f}{\partial
%y}(t,x)+\overline{\nabla}_{\varphi}^{\ast}H(t,x),\\ 
%\]
%\[
%c\left( t,x\right)=\tfrac{\partial \sigma}{\partial
%y}(t,x)
%\]
%Assumptions:
%\begin{description}
%\item[(i)] $\overline{\beta }^{i},\overline{b},c$ and their derivatives $%
%w,r,t,x$ up to order$\ m$ and $\overline{\alpha }^{ij}$ up to order $\max
%\left\{ 2,m\right\} $ are uniformly bounded by the positive constant $k.$

%\item[(ii)] For each $\left( \omega,t,x\right) \in \Omega \times \left[ 0,T\right]
%\times \mathbb{R}$%
%\[
%2\overline{\alpha }^{ij}\left( r,x\right) \xi ^{i}\xi ^{j}\geq 0,\forall \xi
%\in \mathbb{R} 
%\]
%\end{description}

\subsection{A sufficient maximum principle approach (I)}

We now formulate a sufficient version ( a verification theorem) of the maximum
principle for the optimal control of the problem \eqref{eq2.1}-\eqref{eq2.4}.

\begin{theorem}
[Sufficient Maximum Principle (I)]Suppose $\widehat{u}\in\mathcal{U}$, with
corresponding \newline$\widehat{Y}(t,x),\widehat{p}(t,x),\widehat
{q}(t,x).$ Suppose the functions $(y,\varphi)\mapsto
g(x,y,\varphi)$ and\newline$(y,\varphi,u)\mapsto H(t,x,y,\varphi,u,\widehat
{p}(t,x),\widehat{q}(t,x))$ are convex for each
$(t,x)\in\lbrack0,T]\times D$. Moreover, suppose that, for all $(t,x)\in
\lbrack0,T]\times D$,
\begin{align*}
&  \min_{v\in U}H(t,x,\widehat{Y}(t,x),\widehat{Y}(t,\cdot),v,\widehat
{p}(t,x),\widehat{q}(t,x))\\
&  =H(t,x,\widehat{Y}(t,x),\widehat{Y}(t,\cdot),\widehat{u}(t,x),\widehat
{p}(t,x),\widehat{q}(t,x)).
\end{align*}

Then $\widehat{u}$ is an optimal control.
\end{theorem}

\dproof Consider
\[
J(u)-J(\widehat{u})=I_{1}+I_{2},
\]
where
\[
I_{1}=\mathbb{E}\left[  \int_{0}^{T}{\int_{D}}\{f(t,x,Y(t,x),Y(t,\cdot
),u(t,x))-f(t,x,\widehat{Y}(t,x),\widehat{Y}(t,\cdot),\widehat{u}%
(t,x))\}dxdt\right]  ,
\]
and
\[
I_{2}={\int_{D}}\mathbb{E}\left[  g(x,Y(T,x),Y(T,\cdot))-g(x,\widehat
{Y}(T,x),\widehat{Y}(T,\cdot))\right]  dx.
\]
By convexity on $g$ together with the identities
\eqref{eq2.6a}-\eqref{eq2.8a}, we get
\begin{align*}
I_{2} &  \geq{\int_{D}}\mathbb{E}\left[  \tfrac{\partial\widehat{g}}{\partial
y}(T,x)(Y(T,x)-\widehat{Y}(T,x))+\left\langle \nabla_{\varphi}\widehat
{g}(T,x),(Y(T,\cdot)-\widehat{Y}(T,\cdot))\right\rangle \right]  dx\\
&  ={\int_{D}}\mathbb{E}\left[  \tfrac{\partial\widehat{g}}{\partial
y}(T,x)(Y(T,x)-\widehat{Y}(T,x))+\overline{\nabla} \widehat
{g}(T,x)(Y(T,x)-\widehat{Y}(T,x))\right]  dx\\
&  ={\int_{D}}\mathbb{E}\left[  \widehat{p}(T,x)(Y(T,x)-\widehat
{Y}(T,x))\right]  dx\\
&  ={\int_{D}}\mathbb{E}\left[  \widehat{p}(T,x)\widetilde{Y}(T,x)\right]  dx,
\end{align*}
where we put
\begin{align}
\widetilde{Y} (t,x)=Y(t,x) - \widehat{Y}(t,x); (t,x) \in [0,T] \times D.
\end{align}
Applying the It\^ o formula to $\widehat{p}(t,x)\widetilde{Y}(t,x)$, we have
\begin{align}
I_{2} &  \geq{\int_{0}^{T}}{\int_{D}}\mathbb{E}\Big[  \widehat{p}
(t,x)\{A_{x}\widetilde{Y}(t,x)+\widetilde{b}(t,x)\}-\widetilde{Y}
(t,x)\{A_{x}^{\ast}\widehat{p}(t,x) \nonumber\\
&  +\dfrac{\partial\widehat{H}}{\partial y}(t,x)+\overline{\nabla}_{\varphi
}^{\ast}\widehat{H}(t,x)\}+\widehat{q}(t,x)\widetilde{\sigma}(t,x)\Big]  dxdt,\label{1}
\end{align}
where
\begin{align}
\widetilde{b}(t)=b(t) - \widehat{b}(t), \text{ } \widetilde{\sigma}(t) = \sigma(t) - \widehat{\sigma}(t).
\end{align}
Since $\widetilde{Y}(t,x)=\widehat{p}(t,x)\equiv0$, for all $(t,x)\in(0,T)\times\partial D$, we get
\begin{equation}
\int_{D} \widehat{p}(t,x)A_{x}\widetilde{Y}(t,x) dx=\int_{D} \widetilde
{Y}(t,x)A_{x}^{\ast}\widehat{p}(t,x)dx.\label{**}
\end{equation}

Substituting $\left(  \ref{**}\right)  $ in $\left(  \ref{1}\right)  $, yields%
\begin{align}
I_{2} &  \geq\int_{0}^{T}\int_{D}\mathbb{E}\left[  \widehat{p}(t,x)\widetilde
{b}(t,x)-\widetilde{Y}(t,x)\left\{  \tfrac{\partial\widehat{H}}{\partial
y}(t,x)+\overline{\nabla} \widehat{H}(t,x)\right\} + \widehat{q}(t,x)\widetilde{\sigma}(t,x)\right]
dxdt.\label{I2} 
\end{align}
Using the definition of the Hamiltonian $H$ in  \eqref{h}, and putting
\begin{align}
\widetilde{H}(t,x)&=H(t,x,Y(t,x),Y(t,\cdot),u(t,x),\widehat
{p}(t,x),\widehat{q}(t,x))\nonumber\\
&-H(t,x,\widehat{Y}(t,x),\widehat{Y}(t,\cdot),\widehat{u}(t,x),\widehat
{p}(t,x),\widehat{q}(t,x)),
\end{align}
we get
\begin{align}
I_{1} &  =\mathbb{E}\left[  \int_{0}^{T}\int_{D}\{\widetilde{H}(t,x)-\widehat
{p}(t,x)\widetilde{b}(t,x)-\widehat{q}(t,x)\widetilde{\sigma}(t,x)\}dxdt\right]  \nonumber\\
&  \geq\mathbb{E}\left[  \int_{0}^{T}\int_{D}\left\{  \dfrac{\partial
\widehat{H}}{\partial y}(t,x)\widetilde{Y}(t,x)+\left\langle \nabla \widehat{H}(t,x),\widetilde{Y}(t,\cdot)\right\rangle \right.  \right.
\nonumber\\
& \left.  \left.  +\dfrac{\partial\widehat{H}}{\partial u}(t,x)\widetilde{u}(t,x)-\widehat
{p}(t,x)\widetilde{b}(t,x)-\widehat{q}(t,x)\widetilde{\sigma}(t,x)  \right\}  dxdt\right],\label{I4}%
\end{align}
where the last inequality holds because of the concavity assumption of
$H$.\newline Summing $\left(  \ref{I2}\right)  $ and $\left(  \ref{I4}\right)
$, and using \eqref{eq2.6a}, \eqref{eq2.8a}, we end up with%
\begin{equation*}%
\begin{array}
[c]{ll}%
I_{1}+I_{2} & \geq\mathbb{E}\left[
%TCIMACRO{\dint _{0}^{T}}%
%BeginExpansion
{\displaystyle\int_{0}^{T}}
%EndExpansion%
%TCIMACRO{\dint _{D}}%
%BeginExpansion
{\displaystyle\int_{D}}
%EndExpansion
\dfrac{\partial\widehat{H}}{\partial u}(t,x)\widetilde{u}(t,x)dxdt\right].
\end{array}
\end{equation*}
By the maximum condition of $H$ we have%
\[
J(u)-J(\widehat{u})\geq\mathbb{E}\left[  \int_{0}^{T}\int_{D}\tfrac
{\partial\widehat{H}}{\partial u}(t,x)\widetilde{u}(t,x)dxdt\right]  \geq0.
\]
\fproof

\subsection{A necessary maximum principle approach (I)}

We now go to the other version of the necessary maximum principle which can be
seen as an extension of Pontryagin's maximum principle to SPDE with space-mean
dynamics. Here concavity assumptions are not required . We consider the
following: \newline Given arbitrary controls $u,\widehat{u}\in\mathcal{U}$
with $u$ bounded, we define
\[%
\begin{array}
[c]{ll}%
u^{\theta}:=\widehat{u}+\theta u; & \theta\in\left[  0,1\right]  .
\end{array}
\]
Note that, thanks to the convexity of $U$, we also have $u^{\theta}%
\in\mathcal{U}$. We denote by $Y^{\theta}:=Y^{u^{\theta}}$ and by $\widehat
{Y}:=Y^{\widehat{u}}$ the solution processes of \eqref{eq2.1} corresponding to
$u^{\theta}$ and $\widehat{u},$\ respectively.\newline\newline Define the
derivative process $Z(t,x)$ by 
\begin{align}
Z(t,x)=\lim_{\theta \rightarrow 0} \frac{1}{\theta}(Y^{\theta}(t,x)-\widehat{Y}(t,x)) \text{ (limit in } \mathcal{Y}^{(T)} ).
\end{align}

Then, by our assumptions on $f,g,b$ and $\sigma$ it is easy to see that $Z(t,x)$ exists and satisfies the following equation:
\small
\begin{equation}
\left\{
\begin{array}
[c]{ll}%
dZ(t,x) & =\left\{  A_{x}Z(t,x)+\dfrac{\partial b}{\partial y}%
(t,x)Z(t,x)+\left\langle \nabla b(t,x),Z(t,\cdot)\right\rangle
+\dfrac{\partial b}{\partial u}(t,x)u(t,x)\right\}  dt\\
& +\left\{  \dfrac{\partial\sigma}{\partial y}(t,x)Z(t,x)+\left\langle
\nabla \sigma(t,x),Z(t,\cdot)\right\rangle +\dfrac{\partial\sigma
}{\partial u}(t,x)u(t,x)\right\}  dB(t),\\
Z(t,x) & =0; \quad(t,x) \in(0,T) \times\partial D,\\
Z(0,x) & =0;\quad x\in D.
\end{array}
\right.  \label{var1}%
\end{equation}
Note that \eqref{var1}, regarded as an SPDE in the unknown $\mathcal{Y}^{(T)}$-valued process $Z$, is linear and hence the existence and uniqueness of solution follows from e.g. Theorem 3.3 in \cite{GM}.

\begin{theorem}
[Necessary Maximum Principle (I)]Let $\widehat{u}(t,x)$ be an optimal control
and $\widehat{Y}(t,x)$ the corresponding trajectory and adjoint processes
$(\widehat{p}(t,x),\widehat{q}(t,x))$. Then we have
\[
\left.  \dfrac{\partial\widehat{H}}{\partial u}\right\vert _{u=\hat{u}%
}(t,x)=0;\quad\text{ a.s. }
\]

\end{theorem}

\dproof Since $\widehat{u}$ is optimal we get, by the definition \eqref{jcost} of $J$, dominated convergence and the chain rule,
\begin{align*}
0&\leq\underset{\theta\rightarrow 0}{\underline{\lim}}\frac{J(u^{\theta
})-J(\widehat{u})}{\theta}\nonumber\\
&= \underset{\theta\rightarrow 0}{\underline{\lim}}\frac{1}{\theta} \mathbb{E}\Big[ \int_{D} \{g(x,Y^{\theta}(T),Y^{\theta}(T,\cdot)) - g(x,\widehat{Y}(T,x),\widehat{Y}(T,\cdot))\} dx \nonumber\\
&+ \int_{D} \int_0^T \{f(t,x,Y^{\theta}(t,x),Y^{\theta}(t, \cdot),u(t,x))-f(t,x,\widehat{Y}(t,x),\widehat{Y}(t, \cdot),u(t,x)) \} dt dx\Big]\nonumber\\
&=  \mathbb{E}\Big[ \int_{D} \underset{\theta\rightarrow 0}{\underline{\lim}}\frac{1}{\theta}\{g(x,Y^{\theta}(T),Y^{\theta}(T,\cdot)) - g(x,\widehat{Y}(T,x),\widehat{Y}(T,\cdot))\} dx \nonumber\\
&+ \int_{D} \int_0^T\underset{\theta\rightarrow 0}{\underline{\lim}}\frac{1}{\theta} \{f(t,x,Y^{\theta}(t,x),Y^{\theta}(t, \cdot),u(t,x))-f(t,x,\widehat{Y}(t,x),\widehat{Y}(t, \cdot),u(t,x)) \} dt dx\Big]\nonumber\\
&= \mathbb{E}\Big[ \int_{D} \frac{\partial g}{\partial y}(x, \widehat{Y}^{\theta}(T,x),\widehat{Y}^{\theta}(T,\cdot))\underset{\theta\rightarrow 0}{\underline{\lim}}\frac{1}{\theta}(Y^{\theta}(t,x)-\widehat{Y}(t,x)) \\
&+\left\langle \nabla g(x, \widehat{Y}^{\theta}(T,x),\widehat{Y}^{\theta}(T,\cdot)),\underset{\theta\rightarrow 0}{\underline{\lim}}\frac{1}{\theta}(Y^{\theta}(t,\cdot)-\widehat{Y}(t,\cdot)) \right\rangle dx \\
&+ \int_{D} \int_0^T \underset{\theta\rightarrow 0}{\underline{\lim}}\frac{1}{\theta}\{f(t,x,Y^{\theta}(t,x),Y^{\theta}(t, \cdot),u(t,x))-f(t,x,\widehat{Y}(t,x),\widehat{Y}(t, \cdot),u(t,x)) \} dt dx\Big]\nonumber\\
\end{align*}
Therefore, writing $\frac{\partial \widehat{g}}{\partial y}(T,x)=\frac{\partial g}{\partial y}(x, \widehat{Y}(T,x),\widehat{Y}(T,\cdot))$ and $\frac{\partial \widehat{f}}{\partial y}(t,x)=\frac{\partial f}{\partial y}(t,x,\widehat{Y}(t,x),\widehat{Y}(t,\cdot), \widehat{u}(t,x))$ and similarly with $\nabla \widehat{g}(T,x),\nabla \widehat{f}(t,x)$ we obtain

\begin{equation}%
\begin{array}
[c]{l}%
0\leq\underset{\theta\rightarrow0}{\underline{\lim}}\frac{J(u^{\theta
})-J(\widehat{u})}{\theta}\\
=\mathbb{E}\left[
%TCIMACRO{\dint _{D}}%
%BeginExpansion
{\displaystyle\int_{D}}
%EndExpansion
\Big \{\dfrac{\partial\widehat{g}}{\partial y}(T,x)Z
(T,x)+\left\langle \nabla \widehat{g}(T,x),Z(T,\cdot
)\right\rangle \Big \}dx\right] \\
+\mathbb{E}\left[
%TCIMACRO{\dint _{0}^{T}}%
%BeginExpansion
{\displaystyle\int_{0}^{T}}
%EndExpansion%
%TCIMACRO{\dint _{D}}%
%BeginExpansion
{\displaystyle\int_{D}}
%EndExpansion
\left\{  \dfrac{\partial\widehat{f}}{\partial y}(t,x)Z
(t,x)+\left\langle \nabla \widehat{f}(t,x),Z(t,\cdot
)\right\rangle +\dfrac{\partial\widehat{f}}{\partial u}(t,x)u(t,x)\right\}
dxdt\right]  .
\end{array}
\label{j}%
\end{equation}
By \eqref{eq2.6a} and the BSPDE for $\widehat{p}(t,x),$ we have
\[
\mathbb{E}\left[  \int_{D}\Big\{\dfrac{\partial\widehat{g}}{\partial
y}(T,x)Z(T,x)+\left\langle \nabla \widehat{g}%
(T,x),Z(T,\cdot)\right\rangle \Big\}dx\right]  =\mathbb{E}\left[
\int_{D}\widehat{p}(T,x)Z(T,x)dx\right]  ,
\]
The It\^ o formula applied to the product $\widehat{p}(t,x)\cdot Z(t,x)$, where $\widehat{p}$ and $Z$ are the associated equations \eqref{var1}, \eqref{p}-\eqref{p_T}, respectively, to the optimal control $\widehat{u}$, combined with the definition of $\widehat{H} 
$ in \eqref{h}, leads to
\small

\begin{align*}
&  \mathbb{E}\left[  \int_{D}\widehat{p}(T,x)Z(T,x)dx\right]
=\mathbb{E}\left[  \int_{0}^{T}\int_{D}\Big (\widehat{p}(t,x)dZ(t,x)+Z(t,x)d\widehat{p}(t,x)\right. \\
& \left.   +\int_{0}^{T}\int_{D}\Big \{\widehat{q}(t,x)\left(  \frac{\partial\sigma
}{\partial y}(t,x)Z(t,x)+\left\langle \nabla 
\sigma(t,x),Z(t,\cdot)\right\rangle +\frac{\partial\sigma}{\partial
u}(t,x)u(t,x)\right) \Big \}dtdx\right] \\
&  =\mathbb{E}\Big [\int_{0}^{T}\int_{D}\Big \{\widehat{p}(t,x)\Big(A_{x}%
Z(t,x)+\frac{\partial b}{\partial y}(t,x)Z
(t,x)+\left\langle \nabla b(t,x),Z(t,\cdot)\right\rangle
+\frac{\partial b}{\partial u}(t,x)u(t,x)\Big )\\
&  +Z(t,x)\Big(-A_{x}^{\ast}\widehat{p}(t,x)-\frac{\partial
\widehat{H}}{\partial y}(t,x)-\overline{\nabla} \widehat
{H}(t,x)\Big)\\
&  +\Big \{\widehat{q}(t,x)\left(  \frac{\partial\sigma}{\partial
y}(t,x)Z(t,x)+\left\langle \nabla \sigma(t,x),Z(t,\cdot)\right\rangle +\frac{\partial\sigma}{\partial u}(t,x)u(t,x)\right)\Big \}dtdx\Big ].
\end{align*}
Substituting this in (\ref{j}), we get
\[
0\leq\mathbb{E}\left[  \int_{0}^{T}\int_{D}\frac{\partial\widehat{H}}{\partial
u}(t,x)u(t,x)dxdt\right]  .
\]
In particular, if we apply this to
\[
u(t,x)=\mathbf{1}_{[s,T]}(t)\alpha(x),
\]
where $\alpha(x)$ is bounded and $\mathcal{F}_{s}$-measurable we get
\[
0\geq\mathbb{E}\left[  \int_{s}^{T}\int_{D}\frac{\partial\widehat{H}}{\partial
u}(t,x)\alpha(x)dxdt\right]  .
\]
Since this holds for all such $\alpha$ (positive or negative) and all
$s\in\lbrack0,T]$, we conclude that
\[
0 = \frac{\partial\widehat{H}}{\partial u}(t,x);\quad\text{ for a.a. }t,x.
\]
\fproof
%$\square$

\subsection{Controls which are independent of $x$}

In many situations, for example in connection with partial observation
control, it is of interest to study the case when the controls
$u(t)=u(t,\omega)$ are not allowed to depend on the space variable $x$. Let us denote the set of such controls $u\in
\mathcal{U}$ by $\overline{\mathcal{U}}$. Then the corresponding control
problem is to find $\widehat{u}\in\overline{\mathcal{U}}$ such that
\[
J(\widehat{u})=\inf_{u\in\overline{\mathcal{U}}}J(u).
\]
The equations for $J$, $Y$, $H$ and $p$ are as before, except that we replace
$u(t,x)$ by $u(t)$. We handle this situation
by introducing integration with respect to $dx$ in the Hamiltonian. We state
the corresponding modified theorems without proofs:

\begin{theorem}
[Sufficient Maximum Principle (II)]Suppose  $\widehat{u}\in\overline
{\mathcal{U}}$, with corresponding \newline$\widehat{Y}(t,x),\widehat
{p}(t,x),\widehat{q}(t,x).$ Suppose the functions
$(y,\varphi)\mapsto g(x,y,\varphi)$ and \newline$(y,\varphi,u)\mapsto
H(t,x,y,\varphi,u,\widehat{p}(t,x),\widehat{q}(t,x)$
are convex for each $(t,x)\in\lbrack0,T]\times D$. Moreover, suppose the
following average minimum condition,
\begin{align*}
&  \min_{v\in U}\left\{  \int_{D}H(t,x,\widehat{Y}(t,x),\widehat{Y}%
(t,\cdot),v,\widehat{p}(t,x),\widehat{q}(t,x))dx\right\} \\
&  =\int_{D}H(t,x,\widehat{Y}(t,x),\widehat{Y}(t,\cdot),\widehat
{u}(t),\widehat{p}(t,x),\widehat{q}(t,x))dx.
\end{align*}
Then $\widehat{u}$ is an optimal control.
\end{theorem}

\begin{theorem}
[Necessary Maximum Principle (II)]Let $\widehat{u}(t)$ be an optimal control
and $\widehat{Y}(t,x)$ the corresponding trajectory and adjoint processes
$(\widehat{p}(t,x),\widehat{q}(t,x))$. Then we have
\[
\left.
%TCIMACRO{\dint _{D}}%
%BeginExpansion
{\displaystyle\int_{D}}
%EndExpansion
\dfrac{\partial\widehat{H}}{\partial u}\right\vert _{u=\widehat{u}%
}(t,x)dx=0;\quad\text{ a.s. }dt\times d\mathbb{P}.
\]

\end{theorem}

\section{Application to vaccine optimisation}
Assume that
the \emph{density $Y(t,x)$ of infected individuals in a  population} in a random/noisy environment changes over time $t$ and space point $x$
according to the following space-interaction reaction-diffusion equation
\small
\begin{equation*}%
\begin{cases}
dY(t,x) & =\dfrac{1}{2}\Delta Y(t,x)dt+\Big(\alpha\overline{Y}%
(t,x)-u(t,x)Y(t,x)\Big)dt+\beta Y(t,x)dB(t),\\
Y(0,x) & = \xi(x)\geq 0;\quad x\in D,\\
Y(t,x) & = \eta(t,x) \geq 0;\quad(t,x)\in(0,T)\times\partial D,
\end{cases}
\end{equation*}
where $\alpha,\beta$ are given constants modelling the effect on the growth $dY(t,x)$ of the term $\overline{Y}$ and of the noise, respectively, and $\overline{Y}(t,x)=G(x,Y(t,\cdot)), $ where, as before, $G$ is a space-averaging operator of the form
\begin{equation*}
G(x,\varphi)=\frac{1}{V(K_{r})}\int_{K_{r}}\varphi(x+y)dy;\quad
\varphi\in L^{2}(D), 
\end{equation*}
with $V(\cdot)$ denoting Lebesgue volume and
\[
K_{r}=\{y\in\mathbb{R}^{n};|y|<r\}
\]
is the ball of radius $r>0$ in $\mathbb{R}^{n}$ centered at $0$.\\
%and
%\[
%\Delta=\sum_{i=1}^{n}\frac{\partial^{2}}{\partial x_{i}^{2}}%
%\]
%is the Laplacian. 
By a slight extension of Theorem \ref{theorem3.2} (see Remark \ref{remark3.3}), we know that $Y(t,x)\geq 0$ for all $t,x$.\\
If $u(t,x)$ represents our vaccine effort rate at $(t,x)$, we
define the total expected cost $J(u)$ of the effort by
\begin{equation*}
J(u) =\mathbb{E}\Big[\frac{\rho}{2} \int_D\int_0^Tu(t,x)^2Y(t,x)dtdx + \int_Dh_0(x) Y(T,x)dx\Big],
\end{equation*}
where $\rho > 0$ is a constant, and $h_0(x)> 0$ is a bounded function. Here we may regard the first quadratic  term as the cost of the vaccination effort, with unit price $\rho$, and the second term as the cost of having remaining infection at time $T$.
In this case the Hamiltonian is
\begin{align*}
H(t,x,y,\overline{y},p,q)=(\alpha\overline{y} -uy )p + \beta y q + \frac{\rho}{2}u^2y,
\end{align*}
and the adjoint equation satisfies
\small
\begin{align}\label{eq250}
\begin{cases}
& dp(t,x) =  -\Big[\frac{1}{2} \Delta p(t,x) -u(t,x) p(t,x)+\overline{\nabla}_{\overline{y}%
}H(t,x)+ \beta q(t,x)+ \frac{\rho}{2}u^2(t,x)\Big]dt+ q(t,x) dB(t),\\
& p(T,x) = h_0(x); \quad x \in D \\
& p(t,x) = 0; \quad (t,x) \in (0,T) \times \partial D,
\end{cases}
\end{align}
where, by Example \ref{lem}, $\overline{\nabla}_{\overline{y}} H(t,x)=v_{D}(x) \alpha
p(t,x),$
with $v_{D}(x):=\frac{V((x+K_{r})\cap D)}{V(K_{r})}.$\\
The first order condition for an optimal $u=\widehat{u}$ for $H$ together with the requirement that $Y(t,x)>0.$
lead to 
\[
\widehat{u}(t,x)=\frac{p(t,x)}{\rho}.
\]
Hence the pair of random fields $(\widehat {p},\widehat {q})$ becomes
\small
\begin{equation} 
\begin{cases}
d\widehat {p}(t,x) & =-\Big[  \frac{1}{2}\Delta \widehat {p}(t,x)+ \frac{1}{2\rho} \widehat{p}^2(t,x)+v_{D}(x) \alpha
\widehat {p}(t,x)+\beta \widehat {q}(t,x)\Big]dt+\widehat {q}(t,x)dB(t),\\
\widehat {p}(T,x) & =h_0(x);\quad x\in D,\\
\widehat {p}(t,x) & =0;\quad(t,x)\in(0,T)\times\partial D. 
\end{cases}
\label{eq3.4*}
\end{equation}
Since $h_0$ and all the coefficients of this equation are deterministic, we can conclude that $\widehat{q}=0$ and \eqref{eq3.4*} reduces to the deterministic partial differential equation
\small
\begin{align}
\begin{cases}
\frac{\partial}{\partial t}\widehat {p}(t,x)&= -\Big[ \frac{1}{2}\Delta \widehat {p}(t,x)+ \frac{1}{2\rho} \widehat{p}^2(t,x)+v_{D}(x) \alpha
\widehat {p}(t,x)\Big],\nonumber\\
\widehat {p}(T,x) & =h_0(x);\quad x\in D,\nonumber\\
\widehat {p}(t,x) & =0;\quad(t,x)\in(0,T)\times\partial D.
 \label{3.4a}
 \end{cases}
\end{align}
%\textcolor{red}{
%Define $F(t,x)=v_{D}(x) \alpha
%\widehat {p}(t,x)+\beta \widehat {q}(t,x)$. Then, we observe that $p$ admits the following mild representation
%$$\widehat {p}(t,x) =P_{T-t}\widehat {p}(T,x)+\int_t^T P_{s-t}F(s,x)ds-\int_t^T P_{s-t}\widehat {q}(s,x)dB(s);\quad 0\leq t \leq T,
%$$where $P_t$ denotes the semigroup generated by the Laplacian operator.}

This is a (deterministic) Fujita type backward quadratic reaction diffusion equation. 
We could also from the beginning have allowed $h_0(x)$ to be random and satisfy 
$\mathbb{E}\left[  \int_{D}h_0^{2}(x)dx\right]  <\infty$. Then the equation \eqref{eq3.4*} would have become a nonlinear backward \emph{stochastic} reaction-diffusion equation. We will not discuss this further here, but refer to Bandle, \& Levine \cite{BL}, Dalang et al \cite{DKZ} and Fujita \cite{F} and the references therein for more information. 
\vskip 0.3cm

\thanks{{\bf Acknowledgements}. 
Nacira Agram and Bernt \O ksendal are gratefully acknowledge the financial support provided by the Swedish Research Council grant (2020-04697) and the Norwegian Research Council grant (250768/F20), respectively.}

\end{document}